\pdfoutput=1
\RequirePackage{ifpdf}
\ifpdf 
\documentclass[pdftex]{sigma}
\else
\documentclass{sigma}
\fi

\numberwithin{equation}{section}
\newtheorem{Theorem}{Theorem}[section]
\newtheorem{Corollary}[Theorem]{Corollary}
\newtheorem{Lemma}[Theorem]{Lemma}
\newtheorem{Proposition}[Theorem]{Proposition}
{\theoremstyle{definition}
\newtheorem{Definition}[Theorem]{Definition}
\newtheorem{Remark}[Theorem]{Remark}
}

\DeclareMathOperator\End{End}

\DeclareMathOperator\rank{rank}

\DeclareMathOperator\ad{ad}

\DeclareMathOperator\tr{tr}

\DeclareMathOperator\Res{Res}

\DeclareMathOperator\Span{span}

\begin{document}

\newcommand{\arXivNumber}{1407.5335}

\allowdisplaybreaks

\renewcommand{\thefootnote}{$\star$}

\renewcommand{\PaperNumber}{005}

\FirstPageHeading

\ShortArticleName{Bosonizations of $\widehat{\mathfrak{sl}}_2$ and Integrable Hierarchies}

\ArticleName{Bosonizations of $\boldsymbol{\widehat{\mathfrak{sl}}_2}$ and Integrable Hierarchies\footnote{This paper is
a~contribution to the Special Issue on New Directions in Lie Theory.
The full collection is available at
\href{http://www.emis.de/journals/SIGMA/LieTheory2014.html}{http://www.emis.de/journals/SIGMA/LieTheory2014.html}}}

\Author{Bojko BAKALOV~$^\dag$ and Daniel FLEISHER~$^\ddag$}

\AuthorNameForHeading{B.~Bakalov and D.~Fleisher}

\Address{$^\dag$~Department of Mathematics, North Carolina State University, Raleigh, NC 27695, USA}
\EmailD{\href{mailto:bojko_bakalov@ncsu.edu}{bojko\_bakalov@ncsu.edu}}
\URLaddressD{\url{http://www4.ncsu.edu/~bnbakalo/}}

\Address{$^\ddag$~Faculty of Mathematics and Computer Science, The Weizmann Institute of Science,\\
\hphantom{$^\ddag$}~Rehovot 76100, Israel}
\EmailD{\href{mailto:daniel.fleisher@weizmann.ac.il}{daniel.fleisher@weizmann.ac.il}}

\ArticleDates{Received July 22, 2014, in f\/inal form January 09, 2015; Published online January 14, 2015}

\Abstract{We construct embeddings of $\widehat{\mathfrak{sl}}_2$ in lattice vertex algebras by composing the Wakimoto
realization with the Friedan--Martinec--Shenker bosonization.
The Kac--Wakimoto hierar\-chy then gives rise to two new hierarchies of integrable, non-autonomous, non-linear partial
dif\/ferential equations.
A~new feature of our construction is that it works for any value of the central element of $\widehat{\mathfrak{sl}}_2$;
that is, the level becomes a~parameter in the equations.}

\Keywords{af\/f\/ine Kac--Moody algebra; Casimir element; Friedan--Martinec--Shenker bosonization; lattice vertex algebra;
Virasoro algebra; Wakimoto realization}

\Classification{17B80; 17B69; 37K10; 81R10}

\renewcommand{\thefootnote}{\arabic{footnote}}
\setcounter{footnote}{0}

\vspace{-2mm}

\section{Introduction}
Vertex operators and vertex (operator) algebras are powerful tools for studying inf\/inite-di\-men\-sio\-nal Lie algebras,
their representations, and generalizations~\cite{B,FB,FLM,FK,K2,LL, LW}.
For any even integral lattice~$L$, one constructs the lattice vertex algebra $V_L$ associated to~$L$.
When~$L$ is the root lattice of a~f\/inite-dimensional simply-laced Lie algebra ${\mathfrak{g}}$, this gives the
Frenkel--Kac construction of a~level one representation of the corresponding af\/f\/ine Kac--Moody algebra
$\hat{\mathfrak{g}}$ (see~\cite{FK,K1,K2}).

\looseness=-1
In this paper we construct a~dif\/ferent vertex operator realization of the af\/f\/ine Kac--Moody algebra
$\widehat{\mathfrak{sl}}_2$ of an arbitrary level~$k$.
We start with the Wakimoto realization of $\widehat{\mathfrak{sl}}_2$, which can be viewed as an embedding of the
associated af\/f\/ine vertex algebra in the vertex algebra generated by a~pair of charged free bosons~$a^+$,~$a^-$ (also known
as a~$\beta\gamma$-system) and another free boson~$b$ (which generates the Heisenberg algebra); see~\cite{F, Wak}.
The Friedan--Martinec--Shenker bosonization of~$a^+$,~$a^-$ then gives us an embedding in a~certain lattice vertex algebra
$V_L$ (see~\cite{A, FMS,Wang}).
The resulting realization of $\widehat{\mathfrak{sl}}_2$ has appeared previously in~\cite{FS}, and is also related to
the ones in~\cite{FF,JMX}.
Under some assumptions, we prove the uniqueness of this realization by classifying all such embeddings of
$\widehat{\mathfrak{sl}}_2$ in a~lattice vertex algebra.
We then consider a~twisted representation~$M$ of the vertex algebra~$V_L$ (cf.~\cite{BK,Do,FFR,FLM, Le}), and obtain
on~$M$ a~representation of~$\widehat{\mathfrak{sl}}_2$ of level~$k$.
This representation does not appear to exist elsewhere in the literature and may be of independent interest.

Let~$V$ be a~highest-weight representation of an af\/f\/ine Kac--Moody algebra $\hat{\mathfrak{g}}$, and
$\Omega_2\in\End(V\otimes V)$ be the Casimir operator that commutes with the diagonal action of $\hat{\mathfrak{g}}$
(see~\cite{K1}).
Consider the equation
\begin{gather}
\label{kacwak}
\Omega_2(\tau\otimes\tau)=\lambda\tau\otimes\tau,
\qquad
\tau\in V,
\end{gather}
where \looseness=-1 $\lambda\in\mathbb{C}$ is a~constant such that the equation holds when~$\tau$ is the highest-weight vector $v\in
V$.
Then~\eqref{kacwak} holds for any~$\tau$ in the orbit of~$v$ under the Kac--Moody group associated to
$\hat{\mathfrak{g}}$ (see~\cite{PK}).
Equivalently,~\eqref{kacwak} is satisf\/ied for all~$\tau$ such that $\tau\otimes\tau$ is in the
$\hat{\mathfrak{g}}$-submodule generated by $v\otimes v$.
In the case when~$V$ provides a~vertex operator realization of $\hat{\mathfrak{g}}$, such as the Frenkel--Kac
construction, after a~number of non-trivial changes of variables one can rewrite~\eqref{kacwak} as an inf\/inite sequence
of non-linear partial dif\/ferential equations called the Kac--Wakimoto hierarchy~\cite{KW}.
The action of $\hat{\mathfrak{g}}$ allows one to construct some particularly nice solutions to these equations called
solitons.
For example, the Korteweg--de Vries   and non-linear Schr\"{o}dinger hierarchies are instances of Kac--Wakimoto  
hierarchies related to dif\/ferent realizations of   $\widehat{\mathfrak{sl}}_2$ (see~\cite{K1}).

In this paper we investigate the hierarchy~\eqref{kacwak} arising from the Friedan--Martinec--Shenker bosonization of
the Wakimoto realization of $\widehat{\mathfrak{sl}}_2$, which we call the Wakimoto hierarchy.
The Casimir operator $\Omega_2$ is replaced with one of the operators from the coset Virasoro construction, which still
commutes with the diagonal action of $\widehat{\mathfrak{sl}}_2$ (see~\cite{GKO,KR}).
We write the equations of the Wakimoto hierarchy explicitly as Hirota bilinear equations, and we f\/ind the simplest ones.
This is done both in the untwisted case when~$V$ is a~Fock space contained in $V_L$, and the twisted case when $V=M$ is
a~twisted representation of $V_L$.
The new phenomenon is that these are representations of $\widehat{\mathfrak{sl}}_2$ of any level~$k$, so the level
becomes a~parameter in the equations of the Wakimoto hierarchy.

The paper is organized as follows.
In Section~\ref{swak}, we construct explicitly the embedding of~$\widehat{\mathfrak{sl}}_2$ of level~$k$ in a~lattice
vertex algebra $V_L$, and we prove a~certain uniqueness property of this embedding.
The untwisted Wakimoto hierarchy is investigated in Section~\ref{sint}.
In Section~\ref{stw}, we determine the action of $\widehat{\mathfrak{sl}}_2$ on a~twisted representation~$M$ of $V_L$,
and study the corresponding twisted Wakimoto hierarchy.
The conclusion is given in Section~\ref{sconc}.
Throughout the paper, we work over the f\/ield of complex numbers.

\section{Wakimoto realization and its FMS bosonization}
\label{swak}

We assume the reader is familiar with the basic def\/initions and examples of vertex algebras, and refer to~\cite{FB,FLM, K2,
LL} for more details.
Let us review the \emph{Wakimoto realization} of $\widehat{\mathfrak{sl}}_2$ of level~$k$ from~\cite{Wak} in the
formulation of~\cite{F}.
Consider a~pair of charged free bosons~$a^+$,~$a^-$ (also known as a~$\beta\gamma$-system) and a~free boson~$b$ with the
only non-zero OPEs given by
\begin{gather*}
a^+(z)a^-(w)\sim\frac{1}{z-w},
\qquad
a^-(z)a^+(w)\sim-\frac{1}{z-w}
\qquad \text{and}\qquad
b(z)b(w)\sim\frac{2k}{(z-w)^2}.
\end{gather*}
Then we have a~representation of $\widehat{\mathfrak{sl}}_2$ of level~$k$ def\/ined by
\begin{gather*}
e(z) =a^-(z),
\qquad
h(z) =-2{:}a^+(z)a^-(z){:}+b(z),
\\
f(z) =-2{:}a^+(z)^2a^-(z){:}+k\partial_za^+(z)+{:}a^+(z)b(z){:},
\end{gather*}
where the normal ordering of several terms is from the right.

Recall that the boson-boson correspondence, or \emph{Friedan--Martinec--Shenker $($FMS$)$ bosonization}, provides
a~realization of the charged free bosons~$a^+$,~$a^-$ in terms of f\/ields in a~lattice vertex algebra.
Namely, consider the lattice $\mathbb{Z}\alpha_1+\mathbb{Z}\alpha_2$ with $(\alpha_1|\alpha_2)=0$ and $|\alpha_1|^2 = -
|\alpha_2|^2=1$; then{\samepage
\begin{gather*}
a^+(z) = e^{\alpha_1+\alpha_2}(z),
\qquad
a^-(z) = -{:}\alpha_1(z)e^{-\alpha_1-\alpha_2}(z){:}
\end{gather*}
(see~\cite{A, FMS,Wang}).
Notice that $|\alpha_1+\alpha_2|^2=0$.}

We apply FMS-bosonization to the Wakimoto realization of $\widehat{\mathfrak{sl}}_2$, and we obtain an embedding in
a~lattice vertex algebra $V_L$ of the general form
\begin{gather}
\label{sl2emb}
e(z)={:}\alpha(z)e^\delta(z){:},
\qquad
h(z)= \beta(z),
\qquad
f(z)={:}\gamma(z)e^{-\delta}(z){:},
\end{gather}
where $|\delta|^2=0$ and $\alpha(z)$, $\beta(z)$, $\gamma(z)$ are Heisenberg f\/ields.
We will classify all such embeddings, where the lattice~$L$ contains~$\delta$ and
\begin{gather*}
{\mathfrak{h}}=\mathbb{C}\otimes_\mathbb{Z} L=\Span_\mathbb{C}\{\alpha,\beta,\gamma,\delta\}.
\end{gather*}
Note that the rescaling
\begin{gather}
\label{resc}
e\mapsto\lambda e,
\qquad
f\mapsto\frac1\lambda f,
\qquad
h\mapsto h,
\qquad
\lambda\in\mathbb{C},
\end{gather}
is an automorphism of~${\mathfrak{sl}}_2$.

\begin{Theorem}
\label{twak}
Up to the rescaling~\eqref{resc}, the above formulas~\eqref{sl2emb} provide an embedding of
$\widehat{\mathfrak{sl}}_2$ of level~$k$ in the lattice vertex algebra~$V_L$ if and only if
\begin{gather*}
\beta=k\delta+\alpha-\gamma,
\qquad
|\delta|^2=0,
\qquad
(\alpha|\gamma) = k+1,
\qquad
|\alpha|^2 = |\gamma|^2 = (\delta|\alpha)=-(\delta|\gamma)=1.
\end{gather*}
\end{Theorem}

\begin{Remark}
The matrix of $(\cdot|\cdot)$ relative to the spanning vectors $\{\delta,\alpha,\gamma\}$ of ${\mathfrak{h}}$ is
\begin{gather*}
G=
\begin{pmatrix}
0 & 1 & -1
\\
1 & 1 & k+1
\\
-1 & k+1 & 1
\end{pmatrix}
.
\end{gather*}
Since $\det G=-2k-4$, we have $\dim{\mathfrak{h}}=\rank L=3$ if the level $k\ne-2$ is not critical.
When $k=-2$ is critical, we have $\gamma=-\alpha$ and $\dim{\mathfrak{h}}=\rank L=2$.
\end{Remark}

\begin{proof}[Proof of Theorem~\ref{twak}] We need to verify the following OPEs:
\begin{alignat*}{3}
& e(z)f(w) \sim\frac{h(w)}{z-w}+\frac{k}{(z-w)^2}, \qquad && h(z)h(w) \sim\frac{2k}{(z-w)^2},&
\\
& h(z)e(w) \sim\frac{2e(w)}{z-w},\qquad &&  h(z)f(w) \sim\frac{-2f(w)}{z-w},&
\\
& e(z)e(w) \sim 0, \qquad && f(z)f(w) \sim 0.&
\end{alignat*}
We compute (where ``h.o.t.'' stands for higher order terms in $z-w$):
\begin{gather*}
e(z) f(w)={:}\alpha(z)e^\delta(z){:} \, {:}\gamma(w)e^{-\delta}(w){:}
\\
\hphantom{e(z) f(w)}{} \sim\left(\frac{-(\alpha|\delta)\gamma(w)-(\gamma|\delta)\alpha(z)}{z-w}+\frac{(\alpha|\gamma)+(\delta|\alpha)(\delta|\gamma)}{(z-w)^2}\right)
\bigl(1+(z-w)\delta(w)+\text{h.o.t.}\bigr)
\\
\hphantom{e(z) f(w)}{}
 \sim\frac{\bigl((\alpha|\gamma)+(\delta|\alpha)(\delta|\gamma)\bigr)\delta(w)}{z-w} -
\frac{(\alpha|\delta)\gamma(w)+(\gamma|\delta)\alpha(w)}{z-w}
+\frac{(\alpha|\gamma)+(\delta|\alpha)(\delta|\gamma)}{(z-w)^2}.
\end{gather*}
This implies
\begin{gather*}
k =(\alpha|\gamma)+(\delta|\alpha)(\delta|\gamma),
\qquad
h =\beta=k\delta-(\alpha|\delta)\gamma-(\gamma|\delta)\alpha.
\end{gather*}
Similar computations for $e(z)e(w)$ and $f(z)f(w)$ give
\begin{gather*}
(\alpha|\alpha)=(\delta|\alpha)^2,
\qquad
(\gamma|\gamma)=(\delta|\gamma)^2.
\end{gather*}
Now,
\begin{gather*}
h(z)e(w) =h(z) \, {:}\alpha(w)e^\delta(w){:}
 \sim\frac{(h|\delta){:}\alpha(w)e^\delta(w){:}}{z-w}+\frac{(h|\alpha)}{(z-w)^2}
\end{gather*}
tells us that
\begin{gather*}
(h|\delta)=2,
\qquad
(h|\alpha)=0,
\end{gather*}
and a~nearly identical computation for $h(z)f(w)$ gives
\begin{gather*}
(h|\gamma)=0.
\end{gather*}
Then one checks
\begin{gather*}
(h|h)=(h|k\delta-(\alpha|\delta)\gamma-(\gamma|\delta)\alpha)=k(h|\delta)=2k,
\end{gather*}
as required.
Expanding $(h|\delta)=2$, we f\/ind
\begin{gather*}
2=(h|\delta) =(k\delta-(\delta|\alpha)\gamma-(\delta|\gamma)\alpha|\delta)=-2(\delta|\alpha)(\delta|\gamma),
\end{gather*}
i.e.,
\begin{gather*}
(\delta|\alpha)(\delta|\gamma)=-1.
\end{gather*}
Similarly, from $(h|\alpha)=0$ we get
\begin{gather*}
0 =(h|\alpha) =(k\delta-(\delta|\alpha)\gamma-(\delta|\gamma)\alpha|\alpha)
 = k(\delta|\alpha)-(\delta|\alpha)(\gamma|\alpha)-(\delta|\gamma)(\alpha|\alpha)
\\
\hphantom{0}{} = (\delta|\alpha)\bigl(k+1-(\alpha|\gamma)\bigr),
\end{gather*}
which gives
\begin{gather*}
(\alpha|\gamma)=k+1.
\end{gather*}
Gathering all the lattice equations so far obtained,
\begin{gather*}
(\delta|\alpha)(\delta|\gamma)=-1,
\qquad
(\alpha|\gamma)=k+1,
\qquad
(\alpha|\alpha)=(\delta|\alpha)^2,
\qquad
(\gamma|\gamma)=(\delta|\gamma)^2,
\end{gather*}
we notice that we are free to rescale $\alpha\mapsto \lambda\alpha$ and $\gamma\mapsto\frac{1}{\lambda}\gamma$, which
allows us to f\/ix $(\delta|\alpha)=1$.
This then immediately f\/ixes all the other inner products and we obtain the desired result.
\end{proof}

\begin{Remark}
The above embedding of $\widehat{\mathfrak{sl}}_2$ in $V_L$ is essentially equivalent to the ``symmet\-ric''~$W_2^{(2)}$
algebra of~\cite{FS}, to which they refer as the ``three-boson realization'' of~$\widehat{\mathfrak{sl}}_2$.
\end{Remark}

We will expand f\/ields $\phi(z)$ in the standard way
\begin{gather*}
\phi(z) = \sum\limits_{n\in\mathbb{Z}} \phi_{(n)} z^{-n-1}
\end{gather*}
and call the coef\/f\/icients $\phi_{(n)}$ the modes of $\phi(z)$.
The lattice vertex algebra $V_L$ is equipped with the standard action of the Virasoro algebra (see, e.g.,~\cite{K2}):
\begin{gather*}
L_n = \frac12 \sum\limits_{i=1}^{\rank L} \sum\limits_{m\in\mathbb{Z}} {:}a^i_{(m)} b^i_{(n-m)}{:},
\qquad
n\in\mathbb{Z},
\end{gather*}
where $\{a^i\}$ and $\{b^i\}$ are dual bases of ${\mathfrak{h}}$ with respect to the bilinear form $(\cdot|\cdot)$.

\begin{Lemma}
\label{lvir}
When the level $k=c-1\ne-2$ is not critical, we have
\begin{gather*}
L_n = \frac{c-1}4 \sum\limits_{m\in\mathbb{Z}} {:} \delta_{(m)} \delta_{(n-m)} {:} + \frac12
\sum\limits_{m\in\mathbb{Z}} {:} \delta_{(m)} (\alpha-\gamma)_{(n-m)} {:} \\
\hphantom{L_n =}{}
+ \frac1{4(c+1)} \sum\limits_{m\in\mathbb{Z}}
{:}(\alpha+\gamma)_{(m)} (\alpha+\gamma)_{(n-m)}{:}.
\end{gather*}
This formula remains true for $c=-1$ if we remove the last term and set $\gamma=-\alpha$.
\end{Lemma}
\begin{proof}
The proof is straightforward, using that for $c\ne-1$ the Gram matrix~$G$ of $(\cdot|\cdot)$ relative to the basis
$\{\delta,\alpha,\gamma\}$ of ${\mathfrak{h}}$, and its inverse are
\begin{gather*}
G=
\begin{pmatrix}
0 & 1 & -1
\\
1 & 1 & c
\\
-1 & c & 1
\end{pmatrix}
,
\qquad
G^{-1}=\frac1{2(c+1)}
\begin{pmatrix}
c^2-1 & c+1 & -c-1
\\
c+1 & 1 & 1
\\
-c-1 & 1 & 1
\end{pmatrix}
.
\end{gather*}
For $c=-1$, the pair of dual bases for ${\mathfrak{h}}$ are: $\{\delta,\alpha\}$ and $\{-\delta+\alpha,\delta\}$.
\end{proof}

Note that
\begin{gather}
\label{virf}
Y(\omega,z)=L(z)=\sum\limits_{n\in\mathbb{Z}} L_n z^{-n-2},
\end{gather}
where the \emph{Virasoro vector} $\omega\in V_L$ is given~by
\begin{gather}
\label{om}
\omega = \frac{c-1}4 \delta_{(-1)} \delta + \frac12 \delta_{(-1)} (\alpha-\gamma) + \frac1{4(c+1)}
(\alpha+\gamma)_{(-1)} (\alpha+\gamma).
\end{gather}

\begin{Remark}
The above representation of the Virasoro algebra has a~central charge equal to~$\rank L$ (see, e.g.,~\cite{K2}).
On the other hand, the Sugawara construction provides a~representation of the Virasoro algebra from
$\widehat{\mathfrak{sl}}_2$ of level $k\ne-2$, with central charge equal to $3k/(k+2)$ (see, e.g.,~\cite{K2}).
\end{Remark}

\section{The untwisted Wakimoto hierarchy}
\label{sint}

In the previous section, we saw that the modes of
\begin{gather*}
e(z)  = {:}\alpha(z)e^\delta(z){:},
\qquad
f(z)  = {:}\gamma(z)e^{-\delta}(z){:},
\qquad
h(z)  = k\delta(z)+\alpha(z)-\gamma(z)
\end{gather*}
give a~representation of the af\/f\/ine Kac--Moody algebra $\widehat{\mathfrak{sl}}_2$ of level~$k$, where
\begin{gather*}
|\delta|^2=0,
\qquad
|\alpha|^2 = |\gamma|^2 = (\delta|\alpha)=-(\delta|\gamma)=1,
\qquad
(\alpha|\gamma) = c:= k+1.
\end{gather*}
Introduce the bosonic Fock space
\begin{gather*}
B=\mathbb{C}\big[x,y,t;q,q^{-1}\big],
\end{gather*}
where
\begin{gather*}
x=(x_1,x_2,x_3,\dots),
\qquad
y=(y_1,y_2,y_3,\dots),
\qquad
t=(t_1,t_2,t_3,\dots).
\end{gather*}
The Heisenberg f\/ields $\alpha(z)$, $\gamma(z)$ and $\delta(z)$ act on~$B$ as follows ($n>0$):
\begin{alignat*}{4}
& \alpha_{(n)} =\partial_{x_n}+c \partial_{y_n}+\partial_{t_n}, \qquad && \alpha_{(-n)} =nx_n, \qquad && \alpha_{(0)} =q \partial_{q}, &
\\
& \gamma_{(n)} =c \partial_{x_n}+\partial_{y_n}-\partial_{t_n}, \qquad && \gamma_{(-n)} =ny_n, \qquad && \gamma_{(0)} =-q \partial_{q},&
\\
& \delta_{(n)} =\partial_{x_n}-\partial_{y_n}, \qquad  && \delta_{(-n)} =nt_n,\qquad &&  \delta_{(0)} =0.&
\end{alignat*}
By setting $q=e^\delta$, we identify~$B$ as a~subspace of $V_L$.
Then~$B$ is preserved by the actions of $\widehat{\mathfrak{sl}}_2$ and Virasoro.

Introduce the \emph{Casimir field}, whose modes act on $B\otimes B$ inside $V_L\otimes V_L$,
\begin{gather*}
\Omega(z)=e(z)\otimes f(z)+f(z)\otimes e(z)+\frac{1}{2}h(z)\otimes h(z)-k\otimes L(z)-L(z)\otimes k,
\end{gather*}
where $L(z)$ is the Virasoro f\/ield from the previous section (see~\eqref{virf}).
Note that $\Omega(z)=Y(\Omega,z)$ for
\begin{gather}
\label{Omega}
\Omega = e\otimes f+f\otimes e+\frac12 h\otimes h - k\otimes\omega-\omega\otimes k,
\end{gather}
with~$\omega$ given by~\eqref{om}.

\begin{Proposition}
\label{pcas}
All modes of $\Omega(z)$ commute with the diagonal action of $\widehat{\mathfrak{sl}}_2$, i.e.,
\begin{gather*}
[\Omega(z),a(w)\otimes 1+1\otimes a(w)] = 0,
\qquad
a\in{\mathfrak{sl}}_2.
\end{gather*}
\end{Proposition}
\begin{proof}
This follows from the observation that the modes of $\Omega(z)$ give rise to the coset Virasoro construction
(see~\cite{GKO, KR}).
It can also be derived from the commutator formula from the theory of vertex algebras (see, e.g.,~\cite{K2}):
\begin{gather*}
[A_{(m)},\Omega_{(n)}] = \sum\limits_{j=0}^\infty \binom{m}{j} (A_{(j)}\Omega)_{(m+n-j)},
\end{gather*}
applied to the elements $A=a\otimes1+1\otimes a$ with $a\in{\mathfrak{sl}}_2$.
Using that for $a,b\in{\mathfrak{sl}}_2$ and $j\ge0$,
\begin{gather*}
a_{(j)}\omega = \delta_{j,1} a,
\qquad
a_{(j)}b =
\begin{cases}
[a,b],
&
j=0,
\\
k\tr(ab),
&
j=1,
\\
0,
&
j\ge2,
\end{cases}
\end{gather*}
one easily checks that $A_{(j)}\Omega=0$ for all $j\ge0$.
Therefore, $[A_{(m)},\Omega_{(n)}] = 0$ for $m,n\in\mathbb{Z}$.
\end{proof}

Note that, up to adding a~scalar multiple of the identity operator, $\Omega_{(1)}=\Res_z z\Omega(z)$ is the Casimir
element of $\widehat{\mathfrak{sl}}_2$ (see~\cite{K1}, where it is denoted $\Omega_2$).
The \emph{Kac--Wakimoto hierarchy}~\cite{KW} is given by the equation
$\Omega_{(1)}(\tau\otimes\tau)=\lambda(\tau\otimes\tau)$, where~$\tau$ is in a~certain highest-weight module and
$\lambda\in\mathbb{C}$ is a~constant such that the equation holds when~$\tau$ is the highest-weight vector.
Instead of $\Omega_{(1)}$, we will consider the operator $\Omega_{(0)}=\Res_z \Omega(z)$.
Since the highest-weight vector ${\boldsymbol{1}}\in B$ satisf\/ies
$\Omega_{(0)}({\boldsymbol{1}}\otimes{\boldsymbol{1}})=0$, we obtain the following.

\begin{Corollary}
\label{ccas}
Every vector $\tau\in B$, such that $\tau\otimes\tau$ is in the $\widehat{\mathfrak{sl}}_2$-submodule of
$B\otimes B$ generated by ${\boldsymbol{1}}\otimes{\boldsymbol{1}}$, satisfies the equation
\begin{gather}
\label{untwak}
\Omega_{(0)}(\tau\otimes\tau)=0.
\end{gather} 
\end{Corollary}

We will call~\eqref{untwak} the \emph{untwisted Wakimoto hierarchy}.
Our goal now is to compute explicitly the action of $\Omega_{(0)}$ on $B\otimes B$.
The main step is to simplify $e(z)\otimes f(z)$.
We will use the shorthand notation
\begin{gather*}
\phi' = \phi\otimes 1,
\qquad
\phi'' = 1\otimes\phi,
\end{gather*}
so for example, $x_n' = x_n\otimes1$; and we identify
\begin{gather*}
B\otimes B=\mathbb{C}\big[x',y',t',x'',y'',t'';(q')^{\pm1},(q'')^{\pm1}\big].
\end{gather*}
Recall the elementary Schur polynomials def\/ined by the expansion
\begin{gather}
\label{schur}
\exp\left(\sum\limits_{n=1}^\infty t_n z^n\right)=\sum\limits_{m\in\mathbb{Z}}S_m(t)z^m.
\end{gather}
Clearly, $S_m(t)=0$ for $m<0$, but it will be convenient to use summation over all integers.
Explicitly, one has $S_0(t)=1$ and
\begin{gather*}
S_m(t)=\sum\limits_{m_1+2m_2+3m_3+\cdots=m}\frac{t_1^{m_1}}{m_1!}\frac{t_2^{m_2}}{m_2!}\frac{t_3^{m_3}}{m_3!}\cdots,
\qquad
m\ge1.
\end{gather*}
\begin{Lemma}
\label{lefz}
With the above notation, one has
\begin{gather*}
\Res_z e(z)\otimes f(z) = q' (q'')^{-1} \sum\limits_{i,j,l\in\mathbb{Z}} S_l\bigl(2{\bar t}\bigr) \alpha'_{(i)}
\gamma''_{(j)} S_{l-i-j-1}\bigl(-\tilde{\partial}_{\bar x}+\tilde{\partial}_{\bar y}\bigr),
\end{gather*}
where $2{\bar t}_n = t_n'-t_n''$ and $\tilde{\partial}_{{\bar x}_n} = \frac1n(\partial_{x_n'}-\partial_{x_n''})$.
\end{Lemma}
\begin{proof}
Note that
\begin{gather*}
e^{\pm\delta}(z) =z^{\pm\delta_{(0)}}e^{\pm\delta}\exp\left(\pm\sum\limits_{n>0}\delta_{(-n)}\frac{z^{n}}{n}\right)\exp\left(\pm\sum\limits_{n>0}\delta_{(n)}\frac{z^{-n}}{-n}\right)
\\
\hphantom{e^{\pm\delta}(z) }{}
 =q^{\pm 1}\exp\left(\pm\sum\limits_{n>0}t_n
z^n\right)\exp\left(\pm\sum\limits_{n>0}(\tilde{\partial}_{y_n}-\tilde{\partial}_{x_n})z^{-n}\right),
\end{gather*}
where $\tilde{\partial}_{x_n}=\frac{1}{n}\partial_{x_n}$.
Then
\begin{gather*}
e(z)\otimes f(z) = {:}\alpha'(z) e^{\delta'}(z) \gamma''(z) e^{-\delta''}(z){:} = {:}\alpha'(z) \gamma''(z)
e^{\delta'-\delta''}(z){:}.
\end{gather*}
Expanding the exponentials in $e^{\delta'-\delta''}(z)$, we obtain
\begin{gather*}
e^{\delta'-\delta''}(z)  = q'
(q'')^{-1}\sum\limits_{l,m\in\mathbb{Z}}S_l (t'-t'' )S_m (-\tilde{\partial}_{x'}+\tilde{\partial}_{y'}+\tilde{\partial}_{x''}-\tilde{\partial}_{y''} )z^{l-m}
\\
\hphantom{e^{\delta'-\delta''}(z)}{} = q' (q'')^{-1}\sum\limits_{l,m\in\mathbb{Z}}S_l (2{\bar t} )S_m (-\tilde{\partial}_{\bar
x}+\tilde{\partial}_{\bar y} )z^{l-m},
\end{gather*}
which completes the proof.
\end{proof}

Following the procedure of the Japanese school~\cite{DJKM,DKM, KM} (see also~\cite{KR, MJD}), we will rewrite the
untwisted Wakimoto hierarchy~\eqref{untwak} in terms of Hirota bilinear equations.
Let us recall their def\/inition.

\begin{Definition}
Given a~dif\/ferential operator $P(\partial_x)$ and two functions $f(x)$, $g(x)$, we def\/ine the \emph{Hirota bilinear
operator} $Pf\cdot g$ to be
\begin{gather*}
Pf\cdot g=P(\partial_u) \bigl(f(x+u)g(x-u)\bigr)\big|_{u=0},
\end{gather*}
where $x=(x_1,x_2,\ldots)$, $\partial_x=(\partial_{x_1},\partial_{x_2},\ldots)$, etc.
\end{Definition}

As above, we will consider
\begin{gather*}
\tau\otimes\tau\in\mathbb{C}\big[x',y',t',x'',y'',t'';(q')^{\pm1},(q'')^{\pm1}\big].
\end{gather*}
We introduce the following new variables in $B\otimes B$ (which are dif\/ferent from the previous ones in~$B$):
\begin{alignat*}{3}
& x_n  = \frac12(x_n'+x_n''), \qquad &&  \bar{x}_n  = \frac12(x_n'-x_n''),&
\\
& y_n  = \frac12(y_n'+y_n''), \qquad &&  \bar{y}_n  = \frac12(y_n'-y_n''),&
\\
&t_n  = \frac12(t_n'+t_n''), \qquad &&  \bar{t}_n  = \frac12(t_n'-t_n'').&
\end{alignat*}
Then
\begin{gather*}
 x_n' =x_n+\bar{x}_n,\qquad   x_n'' =x_n-\bar{x}_n,
\qquad
\partial_{x_n'} =\frac{1}{2}(\partial_{x_n}+\partial_{\bar{x}_n}),\qquad
 \partial_{x_n''} =\frac{1}{2}(\partial_{x_n}-\partial_{\bar{x}_n}),
\end{gather*}
and similarly for $y'$, $y''$ and $t'$, $t''$.
Thus
\begin{gather*}
\tau\otimes\tau  = \tau(x',y',t';q') \tau(x'',y'',t'';q'')
 = \tau(x+\bar{x},y+\bar{y},t+\bar{t};q') \tau(x-\bar{x},y-\bar{y},t-\bar{t};q'').
\end{gather*}

In order to rewrite the equations in Hirota bilinear form, we recall the formula (see, e.g.,~\cite{KR}):
\begin{gather}
P(\partial_{\bar x})\tau(x+\bar{x})\tau(x-\bar{x}) =P(\partial_u)\tau(x+\bar{x}+u)\tau(x-\bar{x}-u)|_{u=0}
\nonumber\\
\hphantom{P(\partial_{\bar x})\tau(x+\bar{x})\tau(x-\bar{x})}{}  =Q(\partial_u)\tau(x+u)\tau(x-u)|_{u=0}=Q \tau\cdot\tau,\label{hbo}
\end{gather}
where
\begin{gather*}
Q(\partial_u) = P(\partial_u)\exp\left(\sum\limits_{j=1}^\infty \bar{x}_j\partial_{u_j}\right).
\end{gather*}
Using Lemma~\ref{lefz} and the above notation, we obtain
\begin{gather}
\label{eftau}
\Res_z e(z)\tau\otimes f(z)\tau = q' (q'')^{-1} \sum\limits_{i,j,l\in\mathbb{Z}} S_l\bigl(2{\bar t}\bigr) {:} a'_{(i)}
g''_{(j)} {:}   S_{l-i-j-1}\bigl(-\tilde{\partial}_u+\tilde{\partial}_v\bigr) E   \tau\cdot\tau,
\end{gather}
where
\begin{gather*}
E = \exp\left(\sum\limits_{j=1}^\infty \bar{x}_j\partial_{u_j} + \bar{y}_j \partial_{v_j} + \bar{t}_j \partial_{w_j}
\right),
\\
a'_{(i)} =
\begin{cases}
q'\partial_{q'},
&
i=0,
\\
-i(x_{-i}+\bar{x}_{-i}),
&
i<0,
\\
\frac12(\partial_{x_i}+\partial_{u_i})+\frac{c}2(\partial_{y_i}+\partial_{v_i})+\frac12(\partial_{t_i}+\partial_{w_i}),
&
i>0,
\end{cases}
\end{gather*}
and
\begin{gather*}
g''_{(j)} =
\begin{cases}
-q''\partial_{q''},
&
j=0,
\\
-j(y_{-j}-\bar{y}_{-j}),
&
j<0,
\\
\frac{c}2(\partial_{x_j}-\partial_{u_j})+\frac12(\partial_{y_j}-\partial_{v_j})-\frac12(\partial_{t_j}-\partial_{w_j}),
&
j>0.
\end{cases}
\end{gather*}
Note that $\alpha'_{(i)} = \alpha_{(i)} \otimes 1$ and $\gamma''_{(j)} = 1 \otimes \gamma_{(j)}$ commute, while
$a'_{(i)}$ and $g''_{(j)}$ do not commute for $i=-j$.
As usual, the normally ordered product ${:} a'_{(i)} g''_{(j)} {:}$ is def\/ined by putting all partial derivatives to the
right.

Similarly, by switching the single-primed and double-primed terms, we f\/ind
\begin{gather}
\label{fetau}
\Res_z f(z)\tau\otimes e(z)\tau = (q')^{-1} q'' \sum\limits_{i,j,l\in\mathbb{Z}} S_l\bigl(-2{\bar t}\bigr) {:} a''_{(i)}
g'_{(j)} {:} \, S_{l-i-j-1}\bigl(\tilde{\partial}_u-\tilde{\partial}_v\bigr) E   \tau\cdot\tau,
\end{gather}
where
\begin{gather*}
a''_{(i)} =
\begin{cases}
q''\partial_{q''},
&
i=0,
\\
-i(x_{-i}-\bar{x}_{-i}),
&
i<0,
\\
\frac12(\partial_{x_i}-\partial_{u_i})+\frac{c}2(\partial_{y_i}-\partial_{v_i})+\frac12(\partial_{t_i}-\partial_{w_i}),
&
i>0,
\end{cases}
\end{gather*}
and
\begin{gather*}
g'_{(j)} =
\begin{cases}
-q'\partial_{q'},
&
j=0,
\\
-j(y_{-j}+\bar{y}_{-j}),
&
j<0,
\\
\frac{c}2(\partial_{x_j}+\partial_{u_j})+\frac12(\partial_{y_j}+\partial_{v_j})-\frac12(\partial_{t_j}+\partial_{w_j}),
&
j>0.
\end{cases}
\end{gather*}

The other terms in~\eqref{untwak} are easy to compute.
Recalling that
\begin{gather*}
h_{(j)}=\alpha_{(j)}-\gamma_{(j)}+k\delta_{(j)},
\qquad
k=c-1,
\end{gather*}
we get
\begin{gather*}
\Res_z h(z)\tau \otimes h(z)\tau = \sum\limits_{j\in\mathbb{Z}} h'_{(-j-1)}h''_{(j)}   E   \tau\cdot\tau,
\end{gather*}
where
\begin{gather*}
h'_{(i)} =
\begin{cases}
2q'\partial_{q'},
&
i=0,
\\
-i(x_{-i}+\bar{x}_{-i}) + i(y_{-i}+\bar{y}_{-i}) - ki(t_{-i}+\bar{t}_{-i}),
&
i<0,
\\
\partial_{t_i}+\partial_{w_i},
&
i>0,
\end{cases}
\end{gather*}
and
\begin{gather*}
h''_{(j)} =
\begin{cases}
2q''\partial_{q''},
&
j=0,
\\
-j(x_{-j}-\bar{x}_{-j}) + j(y_{-j}-\bar{y}_{-j}) - kj(t_{-j}-\bar{t}_{-j}),
&
j<0,
\\
\partial_{t_j}-\partial_{w_j},
&
j>0.
\end{cases}
\end{gather*}

Finally, we observe that $\Res_z L(z)=L_{-1}$ and apply Lemma~\ref{lvir} to f\/ind
\begin{gather*}
\Res_z L(z)   \tau\otimes \tau = \frac{c-1}4   \sum\limits_{j\in\mathbb{Z}} d'_{(j)} d'_{(-1-j)}   E   \tau\cdot\tau
+ \frac12 \sum\limits_{j\in\mathbb{Z}} d'_{(j)} (a'_{(-1-j)}-g'_{(-1-j)})  E   \tau\cdot\tau
\\
\hphantom{\Res_z L(z)   \tau\otimes \tau =}{} + \frac1{4(c+1)} \sum\limits_{j\in\mathbb{Z}} (a'_{(j)}+g'_{(j)}) (a'_{(-1-j)}+g'_{(-1-j)})  E  \tau\cdot\tau,
\end{gather*}
where
\begin{gather*}
d'_{(j)} =
\begin{cases}
0,
&
j=0,
\\
-j(t_{-j}+\bar{t}_{-j}),
&
j<0,
\\
\frac12(\partial_{x_j}+\partial_{u_j})-\frac12(\partial_{y_j}+\partial_{v_j}),
&
j>0.
\end{cases}
\end{gather*}
Then $\Res_z \tau\otimes L(z)\tau$ is obtained by switching all single-primed terms with double-primed terms, where
\begin{gather*}
d''_{(j)} =
\begin{cases}
0,
&
j=0,
\\
-j(t_{-j}-\bar{t}_{-j}),
&
j<0,
\\
\frac12(\partial_{x_j}-\partial_{u_j})-\frac12(\partial_{y_j}-\partial_{v_j}),
&
j>0.
\end{cases}
\end{gather*}
In this way, we have rewritten all terms from~\eqref{untwak} as Hirota bilinear operators.
We expand~$\tau$~as
\begin{gather*}
\tau=\sum\limits_{m\in\mathbb{Z}}\tau_m(x,y,t) q^m,
\qquad
\tau_m\in\mathbb{C}[x,y,t];
\end{gather*}
then
\begin{gather}
\label{tautau}
\tau\otimes\tau=\sum\limits_{m,n\in\mathbb{Z}} \tau_m(x',y',t') \tau_n(x'',y'',t'') \, (q')^m (q'')^n.
\end{gather}
Since the functions $\tau_m$ do not depend on any of the variables $\bar{x}_i$, $\bar{y}_i$, $\bar{t}_i$, $q'$ and
$q''$, all coef\/f\/icients in front of monomials in these variables give Hirota bilinear equations for~$\tau_m$.
Observe that, in order to get $(q')^m (q'')^n$ in~\eqref{untwak}, we need to apply $\Res_z e(z)\otimes f(z)$ to the
summand
\begin{gather*}
\tau_{m-1}(x',y',t') \tau_{n+1}(x'',y'',t'')   (q')^{m-1} (q'')^{n+1}
\end{gather*}
from~\eqref{tautau}.
Similarly, $\Res_z f(z)\otimes e(z)$ is applied to
\begin{gather*}
\tau_{m+1}(x',y',t') \tau_{n-1}(x'',y'',t'')   (q')^{m+1} (q'')^{n-1}
\end{gather*}
On the other hand, $\frac12\Res_z h(z)\otimes h(z)$ and $L_{-1}\otimes 1+1\otimes L_{-1}$ have to be applied to
\begin{gather*}
\tau_{m}(x',y',t') \tau_{n}(x'',y'',t'')   (q')^{m} (q'')^{n}.
\end{gather*}

If we further specialize to $m=n=0$, we get equations for the three functions~$\tau_{-1}$,~$\tau_0$ and~$\tau_1$.
We will f\/ind the simplest such equations after setting
\begin{gather*}
x_i=y_i=t_i=\bar x_i=\bar y_i=\bar t_{i-1}=0,
\qquad
i\ge2.
\end{gather*}
Then~\eqref{eftau} reduces to
\begin{gather*}
\Res_z e(z)\tau_{-1}\otimes f(z)\tau_{1} = \sum\limits_{\substack{i,j\ge-1
\\
i+j\le -1}} {:} a'_{(i)} g''_{(j)} {:} \, S_{-1-i-j}\bigl(-\tilde{\partial}_u+\tilde{\partial}_v\bigr) E_1
\tau_{-1}\cdot\tau_{1},
\end{gather*}
where now $a'_{(0)} = g''_{(0)} = -1$ and
\begin{gather*}
E_1 = \exp\bigl(\bar{x}_1\partial_{u_1} + \bar{y}_1\partial_{v_1} \bigr).
\end{gather*}
Similarly, from~\eqref{fetau} we have
\begin{gather*}
\Res_z f(z)\tau_{1}\otimes e(z)\tau_{-1} = \sum\limits_{\substack{i,j\ge-1
\\
i+j\le -1}} {:} a''_{(i)} g'_{(j)} {:} \, S_{-1-i-j}\bigl(\tilde{\partial}_u-\tilde{\partial}_v\bigr) E_1
\tau_{1}\cdot\tau_{-1},
\end{gather*}
where $a''_{(0)} = g'_{(0)} = -1$.
The remaining terms from~\eqref{untwak} become zero when applied to $\tau_0\otimes\tau_0$.

Note that for any polynomial~$P$, we have
\begin{gather*}
P(\partial_u,\partial_v,\partial_w) \tau_{1}\cdot\tau_{-1} = P(-\partial_u,-\partial_v,-\partial_w)
\tau_{-1}\cdot\tau_{1}.
\end{gather*}
Using this, from the coef\/f\/icient of $1$ in~\eqref{untwak}, we f\/ind
\begin{gather*}
\bigr[x_1y_1(\partial_{u_1}-\partial_{v_1})+(x_1+y_1)\bigl]\tau_{-1}\cdot\tau_{1}=0.
\end{gather*}
Similarly, the coef\/f\/icient of $\bar{x}_1^2$ in~\eqref{untwak} gives the equation
\begin{gather*}
\big[x_1y_1(\partial_{u_1}-\partial_{v_1})\partial_{u_1}^2 + (x_1+y_1)\partial_{u_1}^2 +
2y_1(\partial_{u_1}-\partial_{v_1})\partial_{u_1}+2\partial_{u_1}\big]\tau_{-1}\cdot\tau_{1}=0.
\end{gather*}

\section{The twisted Wakimoto hierarchy}
\label{stw}

We now investigate the integrable hierarchy arising from a~twisted representation of $\widehat{\mathfrak{sl}}_2$.
Recall from Theorem~\ref{twak} the embedding of $\widehat{\mathfrak{sl}}_2$ of level~$k$ in the lattice vertex algebra
$V_L$, where~$L$ is a~lattice containing an element~$\delta$ such that
\begin{gather*}
{\mathfrak{h}}=\mathbb{C}\otimes_\mathbb{Z} L=\Span\{\alpha,\gamma,\delta\}
\end{gather*}
and
\begin{gather*}
|\delta|^2=0,
\qquad
|\alpha|^2 = |\gamma|^2 = (\delta|\alpha)=-(\delta|\gamma)=1,
\qquad
(\alpha|\gamma) = c:= k+1.
\end{gather*}
Explicitly, $\widehat{\mathfrak{sl}}_2$ is realized in $V_L$ as
\begin{gather}
\label{efhvl}
e=\alpha_{(-1)}e^\delta,
\qquad
f=\gamma_{(-1)}e^{-\delta},
\qquad
h=k\delta+\alpha-\gamma.
\end{gather}

Observe that ${\mathfrak{h}}$ has an order $2$ isometry~$\sigma$ given~by
\begin{gather*}
\sigma(\alpha)=\gamma,
\qquad
\sigma(\gamma)=\alpha,
\qquad
\sigma(\delta)=-\delta,
\end{gather*}
which preserves the $\widehat{\mathfrak{sl}}_2$ subalgebra described above.
In fact,~$\sigma$ operates on ${\mathfrak{sl}}_2$ as the involution $\sigma=\exp(\frac{\pi\mathrm{i}}2 \ad_{e+f})$,
which acts as
\begin{gather}
\label{carinv}
\sigma(e)=f,
\qquad
\sigma(f)=e,
\qquad
\sigma(h)=-h.
\end{gather}
Composing the embedding $\widehat{\mathfrak{sl}}_2\hookrightarrow V_L$ with any~$\sigma$-twisted representation of
$V_L$, we will obtain a~representation of $\widehat{\mathfrak{sl}}_2$.

We refer the reader to~\cite{BK,Do,FFR,FLM, Le} for twisted modules over lattice vertex algebras.
Here we will only need a~special case.
Recall f\/irst the \emph{$\sigma$-twisted Heisenberg algebra} $\widehat{\mathfrak{h}}_\sigma$, spanned over $\mathbb{C}$
by a~central element~$I$ and elements $a_{(m)}$ ($a\in{\mathfrak{h}}$, $m\in\frac12\mathbb{Z}$) such that $\sigma a~=
e^{-2\pi\mathrm{i} m} a$ (see, e.g.,~\cite{FLM, KP,Le}).
The Lie bracket on $\widehat{\mathfrak{h}}_\sigma$ is given~by
\begin{gather*}
[a_{(m)},b_{(n)}] = m \delta_{m,-n} (a|b) I,
\qquad
a,b\in{\mathfrak{h}},
\qquad
m,n\in \tfrac12\mathbb{Z}.
\end{gather*}
Let $\widehat{\mathfrak{h}}_\sigma^\ge$ (respectively, $\widehat{\mathfrak{h}}_\sigma^<$) be the subalgebra of
$\widehat{\mathfrak{h}}_\sigma$ spanned by all elements $a_{(m)}$ with $m\geq0$ (respectively, $m<0$).
Consider the irreducible highest-weight $\widehat{\mathfrak{h}}_\sigma$-module $M=S(\widehat{\mathfrak{h}}_\sigma^<)$,
called the \emph{$\sigma$-twisted Fock space}, on which~$I$ acts as the identity operator and
$\widehat{\mathfrak{h}}_\sigma^\ge$ annihilates the highest-weight vector ${\boldsymbol{1}}\in M$.

We will denote by $a^M_{(j)}$ the linear operator on~$M$ induced by the action of $a_{(j)}
\in\widehat{\mathfrak{h}}_\sigma$, and will write the \emph{twisted fields} as
\begin{gather*}
Y^M(a,z)=\sum\limits_{j\in\frac12\mathbb{Z}} a^M_{(j)} z^{-j-1},
\qquad
a^M_{(j)} \in\End M.
\end{gather*}
One of the main properties of twisted f\/ields is the~$\sigma$-equivariance
\begin{gather}
\label{tweq}
Y^M(\sigma a,z)=Y^M\big(a,e^{2\pi\mathrm{i}}z\big):=\sum\limits_{j\in\frac12\mathbb{Z}} a^M_{(j)} e^{-2\pi\mathrm{i} j}
z^{-j-1}.
\end{gather}
In our case, this means that when $\sigma a=a$ the modes $a^M_{(j)}$ are nonzero only for $j\in\mathbb{Z}$.
On the other hand, if $\sigma a=-a$ the modes $a^M_{(j)}$ are nonzero only for $j\in\frac12+\mathbb{Z}$.
Note that the eigenspaces of~$\sigma$ on ${\mathfrak{h}}$ are spanned by~$\delta$, $\alpha-\gamma$ and by
$\alpha+\gamma$.

We have the inner products
\begin{gather*}
\Big|\frac{\alpha\pm\gamma}{2}\Big|^2=c_\pm:=\frac{1\pm c}{2},
\qquad
\left(\frac{\alpha-\gamma}{2}\Big|\delta\right)=1,
\end{gather*}
and all other inner products are zero.
Then we can identify
\begin{gather*}
M=\mathbb{C}[x,t],
\qquad
\text{where}  
\qquad
x=(x_1,x_2,x_3,\dots),
\qquad
t=(t_1,t_3,t_5,\dots)
\end{gather*}
and the action of $\widehat{\mathfrak{h}}_\sigma$ is given by
\begin{gather*}
\delta^M_{(j)} =
\begin{cases}
\partial_{t_{2j}},& j\in\big(\frac12+\mathbb{Z}\big)_{>0},
\\
-jx_{-2j},& j\in\big(\frac12+\mathbb{Z}\big)_{<0},
\end{cases}
\\
\left(\frac{\alpha-\gamma}{2}\right)^M_{(j)} =
\begin{cases}
\partial_{x_{2j}},& j\in\big(\frac12+\mathbb{Z}\big)_{>0},
\\
-j(t_{-2j}+c_- x_{-2j}),& j\in\big(\frac12+\mathbb{Z}\big)_{<0},
\end{cases}
\\
\left(\frac{\alpha+\gamma}{2}\right)^M_{(j)} =
\begin{cases}
\partial_{x_{2j}},& j\in\mathbb{Z}_{>0},
\\
-jc_+ x_{-2j},& j\in\mathbb{Z}_{\le0}.
\\
\end{cases}
\end{gather*}
From here we obtain
\begin{gather*}
\alpha^M_{(j)}=
\begin{cases}
\partial_{x_{2j}},& j\in\big(\frac{1}{2}\mathbb{Z}\big)_{>0},
\\
-j(t_{-2j}+c_- x_{-2j}),& j\in\big(\frac12+\mathbb{Z}\big)_{<0},
\\
-jc_+ x_{-2j},& j\in\mathbb{Z}_{\le0}.
\end{cases}
\end{gather*}

Recall that we also have twisted f\/ields corresponding to $e^{\pm\delta}$ (see~\cite{FLM, Le}):
\begin{gather*}
Y^M\big(e^{\pm\delta},z\big)  = \exp\left(\mp\sum\limits_{j\in\left(\frac12+\mathbb{Z}\right)_{<0}} \delta^M_{(j)}
\frac{z^{-j}}j \right) \exp\left(\mp\sum\limits_{j\in\left(\frac12+\mathbb{Z}\right)_{>0}} \delta^M_{(j)}
\frac{z^{-j}}j \right)
\\
\hphantom{Y^M\big(e^{\pm\delta},z\big)}{}
 = \exp\left(\pm\sum\limits_{j\in\left(\frac12+\mathbb{Z}\right)_{>0}} x_{2j} z^j \right)
\exp\left(\mp\sum\limits_{j\in\left(\frac12+\mathbb{Z}\right)_{>0}} \partial_{t_{2j}}  \frac{z^{-j}}j \right).
\end{gather*}
By def\/inition, these also satisfy the~$\sigma$-equivariance~\eqref{tweq}.
Then the embedding~\eqref{efhvl} allows one to extend $Y^M$ to the generators of $\widehat{\mathfrak{sl}}_2$.

\begin{Lemma}
\label{lymez}
We have
\begin{gather*}
Y^M(e,z)= Y^M\big(\alpha_{(-1)}e^\delta,z\big)={:}Y^M(\alpha,z)Y^M\big(e^\delta,z\big){:}-\frac{1}{2z}Y^M\big(e^\delta,z\big).
\end{gather*}
\end{Lemma}
\begin{proof}
Recall that for $a\in{\mathfrak{h}}$, we have
\begin{gather*}
a_{(0)} e^\delta = (a|\delta) e^\delta,
\qquad
a_{(m)} e^\delta = 0,
\qquad
m>0.
\end{gather*}
It follows from~(3.13) in~\cite{BK} that{\samepage
\begin{gather*}
{:}Y^M(a,z) Y^M\big(e^\delta,z\big){:}  = Y^M(a,z)_{(-1)}Y^M\big(e^\delta,z\big)
 =\sum\limits_{m=0}^{\infty}\binom{p}{m}z^{-m}Y^M\big(a_{(m-1)}e^\delta,z\big)
\\
\hphantom{{:}Y^M(a,z) Y^M\big(e^\delta,z\big){:} }{}
 = Y^M\big(a_{(-1)}e^\delta,z\big) + p z^{-1} (a|\delta) Y^M\big(e^\delta,z\big),
\end{gather*}
where $a\in{\mathfrak{h}}$ and $p\in\{0,\frac12\}$ are such that $\sigma a=e^{2\pi\mathrm{i} p}a$.}

Now for $a=\alpha+\gamma$, we have $p=0$ and
\begin{gather*}
{:}Y^M(\alpha+\gamma,z)Y^M\big(e^\delta,z\big){:}=Y^M\big((\alpha+\gamma)_{(-1)}e^\delta,z\big).
\end{gather*}
Similarly, for $a=\alpha-\gamma$, we have $p=\frac12$ and
\begin{gather*}
{:}Y^M(\alpha-\gamma,z)Y^M\big(e^\delta,z\big){:}= Y^M\big((\alpha-\gamma)_{(-1)}e^\delta,z\big)+z^{-1}Y^M\big(e^\delta,z\big),
\end{gather*}
since $(\alpha-\gamma|\delta)=2$.
Adding these two equations gives the desired result.
\end{proof}

The above lemma can also be used to express $Y^M(f,z)$, since
\begin{gather*}
Y^M(f,z)=Y^M\big(e,e^{2\pi\mathrm{i}}z\big)=Y^M\big(e,e^{-2\pi\mathrm{i}}z\big)
\end{gather*}
by~\eqref{carinv},~\eqref{tweq}.
Recall that $V_L$ has a~Virasoro vector~$\omega$ given by~\eqref{om}.
Then we have an action of the Virasoro algebra on the twisted module~$M$, given~by
\begin{gather*}
Y^M(\omega,z)=L^M(z)=\sum\limits_{n\in\mathbb{Z}} L^M_n z^{-n-2}.
\end{gather*}

\begin{Lemma}
\label{lymlz}
We have
\begin{gather*}
L^M_n = \frac{c-1}4 \sum\limits_{j\in\frac12+\mathbb{Z}} {:} \delta^M_{(j)} \delta^M_{(n-j)} {:} + \frac12
\sum\limits_{j\in\frac12+\mathbb{Z}} {:} \delta^M_{(j)} (\alpha-\gamma)^M_{(n-j)} {:} \\
\hphantom{L^M_n =}{} + \frac1{4(c+1)}
\sum\limits_{i\in\mathbb{Z}} {:}(\alpha+\gamma)^M_{(i)} (\alpha+\gamma)^M_{(n-i)}{:} +\frac18 \delta_{n,0}.
\end{gather*}
\end{Lemma}
\begin{proof}
Recall that for $a,b\in{\mathfrak{h}}$, we have
\begin{gather*}
a_{(m)} b = \delta_{m,1} (a|b),
\qquad
m\ge0.
\end{gather*}
Then, as in the proof of Lemma~\ref{lymez} above,
\begin{gather*}
{:}Y^M(a,z) Y^M(b,z){:}  = Y^M(a,z)_{(-1)}Y^M(b,z)
 =\sum\limits_{m=0}^{\infty}\binom{p}{m}z^{-m}Y^M(a_{(m-1)}b,z)
\\
\hphantom{{:}Y^M(a,z) Y^M(b,z){:}}{}
 = Y^M(a_{(-1)}b,z) + \frac{p(p-1)}2 z^{-2} (a|b),
\end{gather*}
where $p\in\{0,\frac12\}$ is such that $\sigma a=e^{2\pi\mathrm{i} p}a$.
Now computing $Y^M(\omega,z)$, the last term in the above equation will be nonzero only when $a=\delta$ and
$b=\alpha-\gamma$, in which case $p=\frac12$ and $(\delta|\alpha-\gamma)=2$.
\end{proof}

The twisted version of the Casimir f\/ield $\Omega(z)$ from the previous section is (cf.~\eqref{Omega}):
\begin{gather*}
\Omega^M(z)  = Y^M(\Omega,z)
 =Y^M(e,z) \otimes Y^M(f,z) + Y^M(f,z) \otimes Y^M(e,z)
\\
\hphantom{\Omega^M(z)  =}{}
+ \frac{1}{2}Y^M(h,z)\otimes Y^M(h,z) - k \otimes L^M(z) - L^M(z) \otimes k.
\end{gather*}
Note that
\begin{gather*}
Y^M(e,z)  \otimes Y^M(f,z) + Y^M(f,z) \otimes Y^M(e,z)
\\
\qquad {} = Y^M(e,z) \otimes Y^M\big(e,e^{2\pi\mathrm{i}}z\big) + Y^M\big(e,e^{-2\pi\mathrm{i}}z\big) \otimes Y^M(e,z).
\end{gather*}
Therefore, when computing the coef\/f\/icients in front of integral powers of~$z$ in $\Omega^M(z)$, we can replace the f\/irst
two terms
\begin{gather*}
Y^M(e,z) \otimes Y^M(f,z) + Y^M(f,z) \otimes Y^M(e,z),
\qquad \text{with}\qquad
2Y^M(e,z) \otimes Y^M\big(e,e^{2\pi\mathrm{i}}z\big).
\end{gather*}

\begin{Theorem}
The modes of the above twisted fields $Y^M(e,z)$, $Y^M(f,z)$ and $Y^M(h,z)$ provide~$M$ with the structure of an
$\widehat{\mathfrak{sl}}_2$-module of level~$k$.
The modes of $\Omega^M(z)$ commute with the diagonal action of $\widehat{\mathfrak{sl}}_2$ on $M \otimes M$, i.e.,
\begin{gather*}
\big[\Omega^M(z), a^M_{(m)} \otimes 1+1 \otimes a^M_{(m)}\big] = 0,
\qquad
a\in{\mathfrak{sl}}_2,
\qquad
m\in\tfrac12\mathbb{Z}.
\end{gather*}
\end{Theorem}

\begin{proof}
We observe that the commutator formula for the modes of twisted f\/ields,
\begin{gather*}
\big[a^M_{(m)},b^M_{(n)}\big] = \sum\limits_{j=0}^\infty \binom{m}{j} (a_{(j)}b)^M_{(m+n-j)},
\end{gather*}
is just like the commutator formula in the vertex algebra itself,
\begin{gather*}
[a_{(m)},b_{(n)}] = \sum\limits_{j=0}^\infty \binom{m}{j} (a_{(j)}b)_{(m+n-j)},
\end{gather*}
provided that~$a$ is an eigenvector of~$\sigma$ (see, e.g.,~\cite{Do,FFR, FLM}).
However, in the modes $a^M_{(m)}$ of twisted f\/ields, the index~$m$ is allowed to be nonintegral.

We already know from Theorem~\ref{twak} that
\begin{gather*}
[a_{(m)},b_{(n)}] = [a,b]_{(m+n)} + m\delta_{m,-n} (a|b) k,
\qquad
a,b\in{\mathfrak{sl}}_2,
\qquad
m,n\in\mathbb{Z},
\end{gather*}
where $(a|b)=\tr(ab)$.
Thus for $a\in\{h,e+f,e-f\}$, we have:
\begin{gather*}
\big[a^M_{(m)},b^M_{(n)}\big] = [a,b]^M_{(m+n)} + m\delta_{m,-n} (a|b) k,
\qquad
b\in{\mathfrak{sl}}_2,
\qquad
n\in\tfrac12\mathbb{Z},
\end{gather*}
where $m\in\mathbb{Z}$ for $a=e+f$ and $m\in\frac12+\mathbb{Z}$ for $a=h,e-f$.
Since~$\sigma$ is an inner automorphism of~${\mathfrak{sl}}_2$, by~\cite[Theorem~8.5]{K1} the above modes $a^M_{(m)}$
give a~representation of the af\/f\/ine Kac--Moody algebra $\widehat{\mathfrak{sl}}_2$.
The statement about $\Omega^M(z)$ follows again from the commutator formula and Proposi\-tion~\ref{pcas}.
\end{proof}

Note that the vector ${\boldsymbol{1}}\in M$ satisf\/ies
\begin{gather*}
a^M_{(j)} {\boldsymbol{1}} = \big(e^{\pm\delta}\big)^M_{(j)} {\boldsymbol{1}} = 0,
\qquad
a\in{\mathfrak{h}},
\qquad
j\in\tfrac12\mathbb{Z}_{\ge0}.
\end{gather*}
By Lemmas~\ref{lefz} and~\ref{lymlz}, this implies
\begin{gather*}
a^M_{(n)} {\boldsymbol{1}} = L^M_n {\boldsymbol{1}} = 0,
\qquad
a\in{\mathfrak{sl}}_2,
\qquad
n\ge 1,
\end{gather*}
while $e^M_{(0)} {\boldsymbol{1}} =-\frac12{\boldsymbol{1}}$ and $L^M_0 {\boldsymbol{1}} =\frac18{\boldsymbol{1}}$.
In particular,
\begin{gather*}
\Omega^M_{(n)} ({\boldsymbol{1}}\otimes{\boldsymbol{1}}) = \Res_z z^n \Omega^M(z)
({\boldsymbol{1}}\otimes{\boldsymbol{1}}) = 0,
\qquad
n\ge 2.
\end{gather*}
Similarly to Corollary~\ref{ccas}, we have the following.

\begin{Corollary}
Every vector $\tau\in M$, such that $\tau\otimes\tau$ is in the $\widehat{\mathfrak{sl}}_2$-submodule of
$M\otimes M$ generated by ${\boldsymbol{1}}\otimes{\boldsymbol{1}}$, satisfies the equation
\begin{gather}
\label{twwak}
\Omega^M_{(2)}(\tau\otimes\tau)=0.
\end{gather}
\end{Corollary}

We will call~\eqref{twwak} the \emph{twisted Wakimoto hierarchy}.
As in Section~\ref{sint}, we will compute explicitly the action of $\Omega^M_{(2)}$ on $M\otimes M$.
We use the same notation as before regarding primed and double-primed objects,
\begin{gather*}
x_n' = x_n\otimes1,
\qquad
x_n'' = 1\otimes x_n,
\qquad
\text{etc.}
\end{gather*}
Slightly abusing the notation, we make the change of variables
\begin{gather*}
 x_n =\frac{1}{2}(x_n'+x_n''),\qquad   \bar{x}_n =\frac{1}{2}(x_n'-x_n''),
\qquad t_n =\frac{1}{2}(t_n'+t_n''), \qquad    \bar{t}_n =\frac{1}{2}(t_n'-t_n'').
\end{gather*}
Then
\begin{gather*}
  x_n' =x_n+\bar{x}_n, \qquad   x_n'' =x_n-\bar{x}_n,
\qquad \partial_{x_n'} =\frac{1}{2} (\partial_{x_n}+\partial_{\bar{x}_n} ),\qquad
 \partial_{x_n''} =\frac{1}{2} (\partial_{x_n}-\partial_{\bar{x}_n} ),
\end{gather*}
and similarly for $t'$, $t''$.
Introduce the ``reduced'' Schur polynomials
\begin{gather*}
R_m(t)=S_m(t_1,0,t_3,0,t_5,0,\dots),
\end{gather*}
where $S_m(t)$ are the elementary Schur polynomials def\/ined by~\eqref{schur}.
Then we can compute the f\/irst term in the expression
\begin{gather*}
\Omega^M_{(2)} = 2\Res_z z^2   Y^M(e,z)\otimes Y^M\big(e,e^{2\pi\mathrm{i}}z\big) \\
\hphantom{\Omega^M_{(2)} =}{} + \tfrac{1}{2} \Res_z z^2   Y^M(h,z)\otimes
Y^M(h,z) - k \otimes L^M_1 - L^M_1 \otimes k.
\end{gather*}

\begin{Lemma}
\label{ltwcas}
We have
\begin{gather*}
\Res_z z^2  \,Y^M(e,z)\otimes Y^M\big(e,e^{2\pi\mathrm{i}}z\big)
 =\sum\limits_{i,j,l\in\mathbb{Z}} (-1)^j R_l(2\bar{x}) \big(\alpha^M_{(i/2)}\big)' (\alpha^M_{(j/2)})''
R_{l-i-j+2}(-2\tilde{\partial}_{\bar t})\\
\qquad{}  - \frac12\sum\limits_{j,l\in\mathbb{Z}} R_l(2\bar{x}) \bigl(\big(\alpha^M_{(j/2)}\big)'
+(-1)^j\big(\alpha^M_{(j/2)}\big)'' \bigr) R_{l-j+2}(-2\tilde{\partial}_{\bar t}) + \frac14\sum\limits_{l\in\mathbb{Z}}
R_l(2\bar{x})R_{l+2}(-2\tilde{\partial}_{\bar t}),
\end{gather*}
where $2{\bar x}_n = x_n'-x_n''$ and $\tilde{\partial}_{{\bar t}_n} = \frac1n(\partial_{t_n'}-\partial_{t_n''})$.
\end{Lemma}
\begin{proof}
Using Lemma~\ref{lymez}, we obtain
\begin{gather*}
Y^M(e,z)  \otimes Y^M\big(e,e^{2\pi\mathrm{i}}z\big) = {:}Y^M(\alpha',z) Y^M\big(\alpha'',e^{2\pi\mathrm{i}}z\big)
Y^M\big(e^{\delta'-\delta''},z\big){:}
\\
\qquad{} - \frac{1}{2z} \, {:}\bigl(Y^M(\alpha',z) + Y^M\big(\alpha'',e^{2\pi\mathrm{i}}z\big) \bigr) Y^M\big(e^{\delta'-\delta''},z\big){:} +
\frac{1}{4z^2}   Y^M\big(e^{\delta'-\delta''},z\big).
\end{gather*}
Then we expand
\begin{gather*}
Y^M \big(e^{\delta'-\delta''},z\big)
 = \exp\left(\sum\limits_{j\in\left(\frac12+\mathbb{Z}\right)_{>0}} (x'_{2j} - x''_{2j}) z^j \right)
\exp\left(\sum\limits_{j\in\left(\frac12+\mathbb{Z}\right)_{>0}} (-\partial'_{t_{2j}} + \partial''_{t_{2j}})
\frac{z^{-j}}j \right)
\\
 = \exp\left(\sum\limits_{j\in\left(\frac12+\mathbb{Z}\right)_{>0}} 2\bar{x}_{2j} z^j \right)
\exp\left(-\sum\limits_{j\in\left(\frac12+\mathbb{Z}\right)_{>0}} 2\tilde{\partial}_{{\bar t}_{2j}} z^{-j} \right)
 = \sum\limits_{l,m\in\mathbb{Z}} R_l(2\bar{x})R_m(-2\tilde{\partial}_{\bar t}) z^{(l-m)/2}.
\end{gather*}
We f\/inish the proof by f\/inding the coef\/f\/icient of $z^{-3}$ in $Y^M(e,z) \otimes Y^M(e,e^{2\pi\mathrm{i}}z)$.
\end{proof}

Now, as in Section~\ref{sint}, we can express the action of $\Omega^M_{(2)}$ on $\tau\otimes\tau$ in terms of Hirota
bilinear operators using formula~\eqref{hbo}.
The recipe is that $\tau\otimes\tau$ gets replaced with $E \tau\cdot\tau$, where
\begin{gather*}
E = \exp\left(\sum\limits_{j=1}^\infty \bar{x}_j\partial_{u_j} + \sum\limits_{i=0}^\infty \bar{t}_{2i+1}
\partial_{w_{2i+1}} \right),
\end{gather*}
and, accordingly, $\partial_{\bar x}$ is replaced with $\partial_u$, while $\partial_{\bar t}$ is replaced with
$\partial_w$.
Then $(\alpha^M_{(j)})'$ becomes
\begin{gather*}
a'_{(j)}=
\begin{cases}
\frac12(\partial_{x_{2j}} + \partial_{u_{2j}}),& j\in\big(\frac{1}{2}\mathbb{Z}\big)_{>0},
\\
-j (t_{-2j}+\bar t_{-2j}) -j c_- (x_{-2j} + \bar x_{-2j}),& j\in\big(\frac12+\mathbb{Z}\big)_{<0},
\\
-jc_+ (x_{-2j} + \bar x_{-2j}),& j\in\mathbb{Z}_{\le0},
\end{cases}
\end{gather*}
while $(\alpha^M_{(j)})''$ becomes
\begin{gather*}
a''_{(j)}=
\begin{cases}
\frac12(\partial_{x_{2j}} - \partial_{u_{2j}}),& j\in\big(\frac{1}{2}\mathbb{Z}\big)_{>0},
\\
-j (t_{-2j}-\bar t_{-2j}) -j c_- (x_{-2j} - \bar x_{-2j}),& j\in\big(\frac12+\mathbb{Z}\big)_{<0},
\\
-jc_+ (x_{-2j} - \bar x_{-2j}),& j\in\mathbb{Z}_{\le0}.
\end{cases}
\end{gather*}
Putting these together, we obtain from Lemma~\ref{ltwcas}
\begin{gather*}
\Res_z z^2  \,Y^M(e,z)\tau\otimes Y^M\big(e,e^{2\pi\mathrm{i}}z\big)\tau
\\
\qquad{} =\sum\limits_{i,j,l\in\mathbb{Z}} (-1)^j R_l(2\bar{x}) \, {:} a'_{(i/2)} a''_{(j/2)} {:} \,
R_{l-i-j+2}(-2\tilde{\partial}_w) E \tau\cdot\tau
\\
\qquad\quad{}
- \frac12\sum\limits_{j,l\in\mathbb{Z}} R_l(2\bar{x}) \bigl(a'_{(j/2)}+(-1)^j a''_{(j/2)} \bigr)
R_{l-j+2}(-2\tilde{\partial}_w) E \tau\cdot\tau
\\
\qquad\quad{}
+ \frac14\sum\limits_{l\in\mathbb{Z}} R_l(2\bar{x})R_{l+2}(-2\tilde{\partial}_w) E \tau\cdot\tau.
\end{gather*}

Recall that, by~\eqref{tweq},
$
\gamma^M_{(j)} = (-1)^{2j} \alpha^M_{(j)}$,
$
j\in\tfrac12\mathbb{Z}$.
So, to compute the other terms in $\Omega^M_{(2)}$, we just need to replace $(\delta^M_{(j)})'$ with
\begin{gather*}
d'_{(j)}=
\begin{cases}
\frac12(\partial_{t_{2j}} + \partial_{w_{2j}}),& j\in\big(\frac{1}{2}+\mathbb{Z}\big)_{>0},
\\
-j (x_{-2j} + \bar x_{-2j}),& j\in\big(\frac12+\mathbb{Z}\big)_{<0},
\end{cases}
\end{gather*}
and $(\delta^M_{(j)})''$ with
\begin{gather*}
d''_{(j)}=
\begin{cases}
\frac12(\partial_{t_{2j}} - \partial_{w_{2j}}),& j\in\big(\frac{1}{2}+\mathbb{Z}\big)_{>0},
\\
-j (x_{-2j} - \bar x_{-2j}),& j\in\big(\frac12+\mathbb{Z}\big)_{<0}.
\end{cases}
\end{gather*}
Then we obtain
\begin{gather*}
\Res_z z^2 \, Y^M(h,z)\tau\otimes Y^M(h,z)\tau =\sum\limits_{j\in\frac12+\mathbb{Z}} h'_{(j)} h''_{(1-j)}
E \tau\cdot\tau,
\end{gather*}
where
\begin{gather*}
h'_{(j)} = k d'_{(j)} + 2 a'_{(j)},
\qquad
h''_{(j)} = k d''_{(j)} + 2 a''_{(j)},
\qquad
j\in\tfrac12+\mathbb{Z}.
\end{gather*}

Finally, we get from Lemma~\ref{lymlz}
\begin{gather*}
L^M_1\tau \otimes \tau = \frac{c-1}4 \sum\limits_{j\in\frac12+\mathbb{Z}} d'_{(j)} d'_{(1-j)} E \tau\cdot\tau \\
\hphantom{L^M_1\tau \otimes \tau =}{}
+
\sum\limits_{j\in\frac12+\mathbb{Z}} d'_{(j)} a'_{(1-j)} E \tau\cdot\tau + \frac1{c+1} \sum\limits_{i\in\mathbb{Z}}
a'_{(i)} a'_{(1-i)} E \tau\cdot\tau.
\end{gather*}
Similarly, $\tau\otimes L^M_1\tau$ is given by the same formula with all primes replaced with double primes.

This completes the rewriting of the twisted Wakimoto hierarchy in terms of Hirota bilinear operators.
Then, as in Section~\ref{sint}, all coef\/f\/icients in front of monomials in the variables $\bar{x}_i$, $\bar{t}_i$ give
Hirota bilinear equations for~$\tau$.
We will f\/ind the simplest such equation by setting
\begin{gather*}
\bar{x}_i=x_{i+2}=\bar{t}_{2i-1}=t_{2i+1}=0,
\qquad
i\ge1.
\end{gather*}
Then $a'_{(j)}=a''_{(j)}=0$ for $j<-1$ and $d'_{(m)}=d''_{(m)}=0$ for $m<-\frac12$.
Thus we obtain
\begin{gather*}
\Omega^M_{(2)} (\tau  \otimes\tau) =2\sum\limits_{\substack{i,j\ge-2
\\
i+j\le2}} (-1)^j \, {:} a'_{(i/2)} a''_{(j/2)} {:} \, R_{2-i-j}(-2\tilde{\partial}_w)   \tau\cdot\tau
\\
\qquad{} - \sum\limits_{j=\pm1,\pm2} \bigl(a'_{(j/2)}+(-1)^j a''_{(j/2)} \bigr) R_{2-j}(-2\tilde{\partial}_w)  \tau\cdot\tau +
\frac12R_{2}(-2\tilde{\partial}_w)   \tau\cdot\tau
\\
\qquad{} + \frac12 \sum\limits_{j=\pm\frac12,\frac32} h'_{(j)} h''_{(1-j)}  \tau\cdot\tau - \frac{k^2}4
\sum\limits_{j=\pm\frac12,\frac32} \bigl(d'_{(j)} d'_{(1-j)} + d''_{(j)} d''_{(1-j)} \bigr)  \tau\cdot\tau
\\
\qquad{}
- k \sum\limits_{j=\pm\frac12,\frac32} \bigl(d'_{(j)} a'_{(1-j)} + d''_{(j)} a''_{(1-j)} \bigr)  \tau\cdot\tau -
\frac{k}{k+2} \sum\limits_{i=-1,2} \bigl(a'_{(i)} a'_{(1-i)} + a''_{(i)} a''_{(1-i)} \bigr)  \tau\cdot\tau.
\end{gather*}

Note that, after setting all $\bar{x}_i$ and $\bar{t}_i$ equal to zero, we have $a'_{(j)} = a''_{(j)}$ and $d'_{(j)} =
d''_{(j)}$ for $j<0$.
Also $P(\partial_u,\partial_w)\tau\cdot\tau = 0$ for any odd polynomial~$P$, i.e., such that
$P(-\partial_u,-\partial_w)=-P(\partial_u,\partial_w)$.
Let us assume, in addition, that~$\tau$ is independent of $x_3$, $x_4$, $t_3$ and $t_5$.
Then the f\/irst term in the above sum simplif\/ies to
\begin{gather*}
-\frac{8}{45}c_+^2x_2^2\partial^6_{w_1} - \frac{1}{3}(t_1+c_-x_1)^2\partial_{w_1}^4 -
2(t_1+c_-x_1)\partial_{x_1}\partial_{w_1}^2 - 4c_+x_2\partial_{x_2}\partial_{w_1}^2 -
\frac{1}{2}\big(\partial_{x_1}^2-\partial_{u_1}^2\big).
\end{gather*}
The other terms of $\Omega^M_{(2)}$ are easier to compute and add up to
\begin{gather*}
\frac{4}{3}c_+x_2\partial_{w_1}^4  + 2\partial_{u_1}\partial_{w_1} + \partial_{x_2} + \frac{1}{2}\partial_{w_1}^2 +
\frac{k^2}{8}\big(\partial_{t_1}^2-\partial_{w_1}^2\big)
\\
\qquad{} + \frac{k}2(\partial_{t_1}\partial_{x_1}-\partial_{u_1}\partial_{w_1}) + \frac12\big(\partial_{x_1}^2-\partial_{u_1}^2\big)
-\frac{k^2}{8}\big(\partial_{t_1}^2+\partial_{w_1}^2\big)
-\frac{k}{2}(\partial_{t_1}\partial_{x_1}+\partial_{u_1}\partial_{w_1}).
\end{gather*}
Putting these together, we obtain that the coef\/f\/icient in front of $1$ in $\Omega^M_{(2)}$ gives the Hirota bilinear
equation
\begin{gather*}
\left[-\frac{8}{45}c_+^2x_2^2\partial^6_{w_1}  - \frac{1}{3}\bigl((t_1+c_-x_1)^2-4c_+x_2\bigr)\partial_{w_1}^4\right.
\\
\left.\qquad {}- 2(t_1 +c_-x_1)\partial_{x_1}\partial_{w_1}^2 + \frac{1}{4}(2-k^2)\partial_{w_1}^2 + (2-k)\partial_{u_1}\partial_{w_1}
+ \partial_{x_2} \right]\tau\cdot\tau=0.
\end{gather*}

We then employ the change of variables
\begin{gather*}
x_2=t,
\qquad
t_1=x,
\qquad
x_1=y,
\qquad
u=\log(\tau),
\end{gather*}
which allows us to write the above as the evolutionary equation
\begin{gather*}
u_t =    \frac{8}{45}c_+^2t^2\big(u_{xxxxxx}+ 30u_{xxxx}u_{xx}+60u_{xx}^3\big)
+\frac13\big((x+c_-y)^2-4c_+t\big)\big(u_{xxxx}+6u_{xx}^2\big)
\\
\hphantom{u_t =}{}
 +\frac{1}{2}\big(k^2-2\big)u_{xx} +(k-2)u_{xy}+2 (x+c_-y ) (u_{xxy}+2u_{xx}u_y ).
\end{gather*}
Note that at the critical level, $k=-2$, we have $c_+=0$, $c_-=1$, and the above equation becomes
\begin{gather*}
u_t=\frac13(x+y)^2\big(u_{xxxx}+6u_{xx}^2\big) + 2(x+y) (u_{xxy}+2u_{xx}u_y ) + u_{xx} - 4u_{xy}.
\end{gather*}
Another reduction is obtained by assuming $u_y=0$ and letting $y=0$.
In this case, we get an order six non-autonomous non-linear PDE, resembling those found in~\cite{DJKM},
\begin{gather*}
u_t =    \frac{2}{45}(k+2)^2   t^2\big(u_{xxxxxx}+ 30u_{xxxx}u_{xx}+60u_{xx}^3\big)
\\
 \hphantom{u_t =}{} +\frac13\big(x^2-2(k+2)t\big)\big(u_{xxxx}+6u_{xx}^2\big) + \frac{1}{2}\big(k^2-2\big)u_{xx}.
\end{gather*}

\section{Conclusion}\label{sconc}

The Frenkel--Kac construction of the homogeneous realization of the basic representation of~$\widehat{\mathfrak{sl}}_2$
provides an embedding of the af\/f\/ine vertex algebra of~$\widehat{\mathfrak{sl}}_2$ at level~$1$ in the lattice vertex
algebra $V_L$, where~$L$ is the root lattice of ${\mathfrak{sl}}_2$ (see~\cite{FK, K2}).
The other realizations of the basic representation of $\widehat{\mathfrak{sl}}_2$ can be viewed then as twisted modules
over $V_L$ (see~\cite{KP,Le, LW}).
These realizations have spectacular applications, such as combinatorial identities obtained from the Weyl--Kac character
formula and integrable systems obtained from the Casimir element (see~\cite{K1}).
In particular, for the principal realization of~\cite{LW} one obtains the Korteweg--de Vries hierarchy, and for the
homogeneous realization of~\cite{FK} one obtains the non-linear Schr\"odinger hierarchy.

In this paper we have constructed other embeddings of the af\/f\/ine vertex algebra of $\widehat{\mathfrak{sl}}_2$ in
lattice vertex algebras $V_L$, for an arbitrary level~$k$ (now~$L$ is not the root lattice).
Then the twisted modules over $V_L$ provide new vertex operator realizations of $\widehat{\mathfrak{sl}}_2$ at
level~$k$.
It will be interesting to understand their representation theoretic signif\/icance, and to generalize them to other Lie
algebras or superalgebras.
In particular, we hope to do this for $\widehat{{\mathfrak{sl}}}_n$, since in this case the Wakimoto realization is
known explicitly by~\cite{FF1}.
Other relevant works include~\cite{A,FS,FF,JMX,Wang}.

As an application of these new vertex operator realizations of $\widehat{\mathfrak{sl}}_2$, we have obtained two
hierarchies of integrable, non-autonomous, non-linear partial dif\/ferential equations.
A~new feature is that the level~$k$ becomes a~parameter in the equations.
It would be interesting to see if the hierarchies associated to $\widehat{\mathfrak{sl}}_2$ (or more generally
$\widehat{{\mathfrak{sl}}}_n$) are reductions of some larger hierarchy, similarly to how the Gelfand--Dickey hierarchies
are reductions of the KP hierarchy.
A~``bosonic'' analog of the KP hierarchy has been constructed by K.~Liszewski~\cite{Li}, and it might be related to one
of our hierarchies when $k=-1$.
Constructing soliton solutions for our equations is as of yet elusive, and is complicated by the fact that all the
f\/ields are bosonic.

\subsection*{Acknowledgements}
We are grateful to Naihuan Jing and Kailash Misra for many useful discussions, and to the referees for valuable
suggestions that helped us improve the exposition.
The f\/irst author was supported in part by NSA and Simons Foundation grants.

\pdfbookmark[1]{References}{ref}
\LastPageEnding

\end{document}